\documentclass[12pt]{elsarticle}

\usepackage[utf8]{inputenc}

\usepackage{amssymb,amsthm, url,amsmath,graphicx, mathtools}
\usepackage{tikz}
\usepackage{pgfplots}
\usepackage{pgfplotstable}
\usepackage{subcaption}
\usepackage{listings}

\newcommand{\MCM}{\mathcal{M}}

\newcommand{\bfb}{\mathbf{b}}
\newcommand{\bfx}{\mathbf{x}}

\newcommand{\software}[1]{\textsc{#1}}
\newcommand{\pytential}{\software{Pytential}}
\newcommand{\firedrake}{Firedrake}
\newcommand{\petsc}{\software{PETSc}}

\newtheorem{theorem}{Theorem}[section]

\usepackage{lineno}

\usepackage{hyperref}
\begin{document}

\begin{frontmatter}


\title{Exact domain truncation for the Morse-Ingard equations\tnoteref{t1}}

\tnotetext[t1]{This work was supported by NSF SHF-1909176 and SHF-1911019.}
\author[bu]{Robert C. Kirby\corref{cor1}}
\ead{robert\_kirby@baylor.edu}

\author[uiuc]{Xiaoyu Wei}
\ead{xywei@illinois.edu}

\author[uiuc]{Andreas Kl{\"o}ckner}
\ead{andreask@illinois.edu}

\cortext[cor1]{Corresponding author}
\address[bu]{Department of Mathematics, Baylor University; Sid Richardson Science Building; 1410 S.~4th St.; Waco, TX 76706; United States}
\address[uiuc]{Department of Computer Science, University of Illinois at Urbana-Champaign,\\ Champaign, Illinois, United
   States}


\begin{abstract}
  Morse and Ingard~\cite{acoustics} give a coupled system of time-harmonic equations for the temperature and pressure of an excited gas.  These equations form a critical aspect of modeling trace gas sensors.
  Like other wave propagation problems, the computational problem must be closed with suitable far-field boundary conditions.
  Working in a scattered-field formulation, we adapt a nonlocal boundary condition proposed in~\cite{kirby2021finite} for the Helmholtz equation to this coupled system.
  This boundary condition uses a Green's formula for the true solution on the boundary, giving rise to a nonlocal perturbation of standard transmission boundary conditions.
  However, the boundary condition is exact and so Galerkin discretization of the resulting problem converges to the restriction of the exact solution to the computational domain.
  Numerical results demonstrate that accuracy can be obtained on relatively coarse meshes on small computational domains, and the resulting algebraic systems may be solved by GMRES using the local part of the operator as an effective preconditioner.
\end{abstract}


\end{frontmatter}


\section{Introduction}
\label{sec:intro}
Laser absorption spectroscopy is used for detecting trace amounts of gases in diverse application areas such as air quality monitoring, disease diagnosis, and manufacturing~\cite{curl2010quantum,patimisco2014quartz,petersen2017quartz}.
In photoacoustic spectroscopy, a laser is fired between the tines of a small quartz tuning fork, and in the presence of a particular gas, acoustic and thermal waves are generated.  These waves, in turn, interact with the tuning fork to generate an electric signal via pyroelectric and piezoelectric effects.  
Two variants of these sensors are the so-called QEPAS (quartz-enhanced photoacoustic spectroscopy) and ROTADE (resonant optothermoacoustic detection) models~\cite{kosterev2002quartz, kosterev2010resonant}.  
In QEPAS, the acoustic wave dominates the signal, while the thermal wave is more important in ROTADE. 
In many experimental configurations, both effects appear.
While a full model of the sensor is an eventual goal, obtaining efficient and accurate models of the gas itself is a critical step.

Earlier work on modeling this problem~\cite{petraThesis, QEPASdesign, ROTADEdesign} simplified the model to a single PDE that included an empirically-determined damping term to account for otherwise-neglected processes.  
This approach is only accurate in particular regimes, and the empirical corrections depend strongly on geometry as well as physical parameters, which limits the model's utility in an optimal design context.

Hence, work began on coupled models including both thermal and acoustic effects.  
A finite element discretization of the coupled pressure-temperature system of Morse and Ingard~\cite{acoustics} was first addressed in~\cite{pyHPC}, where the difficulty of solving the linear system was noted.  
Kirby and Brennan gave a more rigorous treatment in~\cite{brennan2015finite}, with analysis of the finite element error and preconditioner performance.  Kaderli \emph{et al} derived an analytical solution for the coupled
system in idealized geometry in~\cite{kaderli2017analytic}.  
Their technique involves reformulating the system studied in~\cite{brennan2015finite} by an algebraic
simplification that eliminates the temperature Laplacian from the
pressure equation.  
In~\cite{kirby2020optimal}, this reformulation was seen to lose coercivity but still retain
a G{\aa}rding-type inequality, leading to optimal-order finite element convergence theory and preconditioners.
Work by Safin \emph{et al}~\cite{safin} began a more robust multi-physics study, coupling the Morse-Ingard equations for atmospheric pressure and temperature to heat conduction of the quartz tuning fork, although vibrational effects were still not considered.
They also applied a perfectly-matched layer (PML)~\cite{berenger1994perfectly} to truncate the computational domain, and a Schwarz-type preconditioner that separates out the PML region was used to effectively reduce the cost of solving the linear system.
They also include some favorable comparisons between the computational model and experimental data.  

Previous numerical analysis of this problem in the cited literature has focused on volumetric discretizations based on finite elements.  In~\cite{miskie}, we derived a boundary integral formulation for a scattered-field form of the Morse-Ingard equations.  As with other wave problems, this problem writes the solution as the sum of a Morse-Ingard solution that satisfies the forcing (evaluated by means of a fast volumetric convolution with a Green's function) plus a field that satisfies Morse-Ingard with no volumetric forcing but Neumann data on the tuning fork such that the sum satisfies homogeneous boundary conditions.  We then formulated a second-kind integral equation for the scattered field and approximated it with a boundary integral method.  In this work we return to finite element discretization, but we make use of the results we obtained considering the integral form of the equations to make significant advances in imposing a far-field condition.  

In~\cite{kirby2021finite}, we developed a novel nonlocal boundary condition for  truncating the domain of Helmholtz scattering problems.  
This condition, which uses Green's representation of the solution on the artificial boundary to give a nonlocal Robin-type condition involving layer potentials, is exact -- the solution of resulting BVP agrees exactly with the restriction of the solution of the original problem to the computational domain. 
The variational form of the problem inherits a G{\aa}rding-type inequality from the Helmholtz operator so that a Galerkin finite element method yields optimal asymptotic convergence rates.
Moreover, empirical results suggest that the standard local operator with transmission boundary conditions serves as an excellent preconditioner.
Hence, one only needs to apply the action of the resulting nonlocal operator, say, by a fast multipole method, in order to obtain fast GMRES convergence.

In this paper, extend this approach from the Helmholtz operator to the Morse-Ingard system, making using of several results developed in our boundary-integral formulation in~\cite{miskie}.  
In Section~\ref{sec:MI}, we recall the Morse-Ingard equations.
Then, Section~\ref{sec:BC} addresses far-field boundary conditions for the system and appropriate boundary conditions for domain truncation.  
By means of the transformation to a decoupled Helmholtz system, we are able to state an analogous far-field condition and associated transmission-type condition for the Morse-Ingard system.  
This allows a comparison to the 
\emph{ad hoc} transmission boundary conditions used in~\cite{brennan2015finite, kirby2020optimal}.
Moreover, we can derive an exact analog of the nonlocal Helmholtz boundary condition for Morse-Ingard.  
Although one may directly solve the decoupled Helmholtz equations rather than the coupled form of Morse-Ingard, formulating boundary conditions and directly simulating the coupled system serves several purposes.  First, the pointwise transformation, thought it only involves a $2 \times 2$ matrix, is quite ill-conditioned and seems to limit the accuracy we obtain on fine meshes for the decoupled form compared to the coupled system.
Second, a more complete model of trace gas sensors~\cite{safin2018modeling} involves coupling Morse-Ingard to the tuning fork vibration, which in turn requires modeling the fluid flow.  
Domain truncation will still be required, but coupling of pressure and temperature to the fluid and tuning fork may limit the utility of the decoupled system.  Additionally, as noted in~\cite{miskie}, solving for the acoustic mode while neglecting the thermal mode turns out to be an effective approximation. 
After developing the boundary conditions in Section~\ref{sec:BC}, we derive a finite element formulation for the Morse-Ingard system in Section~\ref{sec:varform}.
We discuss the structure of the linear system and approaches to preconditioning in Section~\ref{sec:precond} and we then provide some numerical in Section~\ref{sec:numres} before offering some final conclusions in Section~\ref{sec:conc}.

\section{The Morse-Ingard equations}
\label{sec:MI}

The Morse-Ingard equations of thermoacoustics are a system of partial differential equations for the temperature and pressure of an excited gas.  The model begins from a time-domain formulation.  After assuming time-periodic forcing and performing nondimensionalization and some algebraic manipulations, we arrive at the form given in~\cite{kaderli2017analytic} and further analyzed in~\cite{kirby2020optimal}:

\begin{equation}
  \label{eq:PDE}
  \begin{split}
  -\MCM \Delta T - i T + i \tfrac{\gamma-1}{\gamma} P & = -S, \\
   \gamma\left( 1 - \tfrac{\Lambda}{\MCM} \right) T
  -\left( 1 - i \gamma \Lambda \right) \Delta P
  - \left[ \gamma \left(1-\tfrac{\Lambda}{\MCM}\right)+
    \tfrac{\Lambda}{\MCM}\right]P
  & = i\gamma\tfrac{\Lambda}{\MCM}S.
  \end{split}
\end{equation}

Here, $T$ and $P$ are the non-dimensional temperature and pressure, respectively, within the gas.  $S$ is a volumetric forcing function, modeling for example a laser pulse. $\gamma$ is the ratio of specific heat of the gas at constant pressure to that at constant volume.  The dimensionless number $\MCM$ measures the ratio of the product of the characteristic thermal conduction scale and forcing frequency to sound speed, and $\Lambda$ does similarly for the viscous length scale.  Typical values of parameters
\begin{equation}
  \begin{split}
    \gamma & =7/5 \\
    \MCM &  =3.664152973215096\cdot 10^{-5} \\
    \Lambda & =5.370572762330994\cdot 10^{-5}
  \end{split}
\end{equation}
are taken as in~\cite{kaderli2017analytic, miskie}.

We let $\Omega^c \subset \mathbb{R}^d$ (with $d=2,3$) be a bounded domain representing the tuning fork, and let its boundary be called $\Gamma$.  The complement of $\overline{\Omega^c}$ will be the domain $\Omega$ on which we pose~\eqref{eq:PDE}.  On $\Gamma$, we impose homogeneous Neumann boundary conditions,
\begin{equation}
\label{eq:BC}
\frac{\partial T}{\partial n} = 0, \ \ \ 
\frac{\partial P}{\partial n} = 0.
\end{equation}
which posits that the tuning fork is thermally insulated from the gas, and that the tuning fork is sound-hard.  More advanced models, in which the gas heats the tuning fork or the acoustic waves couple to tuning fork deformation, generalize this condition~\cite{safin}.

A suitable far-field condition is required to close the model, which requires some appropriate decay at infinity akin to the Sommerfeld radiation condition for the Helmholtz operator.  Numerical methods based on volumetric discretization on a truncated domain have posed either some kind of transmission-type condition~\cite{brennan2015finite,kirby2020optimal} or perfectly-matched layers~\cite{safin}.

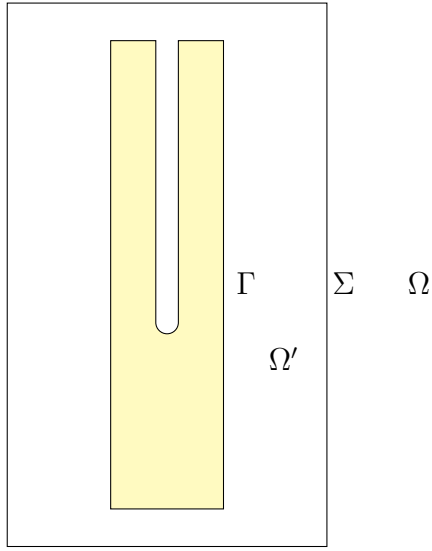
\begin{figure}
\centering
\begin{tikzpicture}[scale=10]
\draw (-0.2125, 0) rectangle (0.2125, 0.723);
\draw[fill=yellow!30] (-0.015, 0.673) -- (-0.075, 0.673) --
(-0.075, 5e-2) -- (0.075, 5e-2) -- (0.075, 0.673) -- (0.015, 0.673) -- (0.015, 0.298)
arc[start angle=0, end angle=-180, radius=0.015] -- cycle;
\node at (0.105, 0.35) {$\Gamma$};
\node at (0.155, 0.25) {$\Omega'$};
\node at (0.235, 0.35) {$\Sigma$};
\node at (0.335, 0.35) {$\Omega$};
\end{tikzpicture}
\caption{
  The yellow-shaded area represents the tuning fork.
  The computational domain $\Omega$ is the rectangle minus the tuning fork.  $\Sigma$ is the outer rectangle, and $\Gamma$ is the boundary of the tuning fork itself.}
\label{fig:tf}
\end{figure}

In~\cite{miskie}, we gave a boundary integral method for Morse-Ingard based on a scattered-field formulation, which turns the volumetric inhomogeneity into an inhomogeneous Neumann condition on~$\Gamma$. We discuss this in greater detail in Section~\ref{sec:morse-ingard-farfield}.
To arrive at the scattered-field formulation, we split the solution into incoming and scattered waves via
\begin{equation}
T = T^i + T^s, \ \ \ P = P^i + P^s,
\end{equation}
where $T^i$ and $P^i$ satisfy~\eqref{eq:PDE} with the given forcing function $S$ but have some inhomogeneous boundary conditions on~$\Gamma$.  Then, $T^s$ and $P^s$ are chosen to satisfy~\eqref{eq:PDE} with homogeneous forcing $S=0$ and such that the combined waves satisfy~\eqref{eq:BC}.
The incoming waves $T^i$ and $P^i$ can be constructed by volumetric convolution of $S$ with a free-space Green's function.  With these in hand, their normal derivatives on $\Gamma$ can be computed, and the negative of these used as boundary conditions for $T^s$ and $P^s$.  Consequently, we drop the superscripts `s' for the scattered field and, for the rest of the paper, consider the system of PDE
\begin{equation}
  \label{eq:scatPDE}
  \begin{split}
  -\MCM \Delta T - i T + i \tfrac{\gamma-1}{\gamma} P & = 0, \\
   \gamma\left( 1 - \tfrac{\Lambda}{\MCM} \right) T
  -\left( 1 - i \gamma \Lambda \right) \Delta P
  - \left[ \gamma \left(1-\tfrac{\Lambda}{\MCM}\right)+
    \tfrac{\Lambda}{\MCM}\right]P
  & = 0,
  \end{split}
\end{equation}
together with boundary conditions
\begin{equation}
\label{eq:scatBC}
\frac{\partial T}{\partial n} = g_T, \ \ \ 
\frac{\partial P}{\partial n} = g_P,
\end{equation}
on $\Gamma$ together with an appropriate far field condition.
Although we obtained satisfactory results with a boundary integral method in~\cite{miskie}, we work with volumetric (finite element) discretizations here to chart a path towards modeling of additional volumetric physical phenomena without additional computational machinery in future work.
Consequently, in this paper we are primarily interested in developing an analog of the nonlocal boundary condition developed in~\cite{kirby2021finite} for the Morse-Ingard system.  To this end, we introduce a truncating boundary $\Sigma$ and define a domain $\Omega' \subset \Omega$ to be that contained between $\Gamma$ and $\Sigma$, as shown in Figure~\ref{fig:tf}.  

In~\cite{miskie}, we demonstrated that the Morse-Ingard system~\eqref{eq:PDE} could be decoupled into a pair of independent Helmholtz equations.  After significant algebraic manipulations, we introduce modified material coefficients
\begin{equation}
\begin{split}
Q^2 & =  4(i\MCM + \gamma \MCM \Lambda) + (1 - i \gamma \MCM - i \Lambda)^2, \\
t_{\pm} & = 
\frac{(2\Lambda\gamma-\Lambda-\MCM\gamma+i)\MCM \mp i\MCM Q }{2\gamma(\Lambda-\MCM)(i\Lambda\gamma-1)}
\end{split}
\end{equation}
as well as separate thermal and acoustic wave numbers
\begin{equation}
\begin{split}
k_t^2 & = 
\frac{i}{2\Omega} \left(\frac{1 - i\gamma\Omega - i\Lambda + Q}{1 -
i\gamma\Lambda}\right), \\
k_p^2 & = \frac{i}{2\Omega} \left(\frac{1 - i\gamma\Omega - i\Lambda - Q}{1 -
i\gamma\Lambda}\right).
\end{split}
\end{equation}
With the change of variables
\begin{equation}
\label{eq:Vdef}
\begin{bmatrix}
V_t \\ V_p
\end{bmatrix}
= \begin{bmatrix} 
\MCM & t_{+} (1-i \gamma \Lambda) \\
\MCM & t_{-} (1-i \gamma \Lambda)
\end{bmatrix}
\begin{bmatrix} T \\ P
\end{bmatrix}
\equiv
B \begin{bmatrix} T \\ P
\end{bmatrix},
\end{equation}
the Morse-Ingard system~\eqref{eq:PDE} decouples into separate
Helmholtz equations
\begin{equation}
\label{eq:MIdecomp}
\begin{split}
\Delta^2 V_t + k_t^2 V_t & = a_t S, \\
\Delta^2 V_p + k_p^2 V_p & = a_p S, \\
\end{split}
\end{equation}
where $a_t$ and $a_p$ are data-dependent constants. For parameter regimes of interest, $k_t$ has a large imaginary part, so that $V_t$ rapidly attenuates.  On the other hand, $k_p$ has a very small imaginary part so that $V_p$ attenuates slowly.  Typical values of these new parameters based on those from given above are
\begin{equation}
  \begin{split}
    k_t & = 116.81449127197266+116.81529235839844i, \\
    k_p & = 1+3.418116830289364e-05i.
  \end{split}
\end{equation}

The same change of variables decouples the scattered-field formulation~\eqref{eq:scatPDE}.  In this case, we have that
\begin{equation}
\label{eq:scatMIdecomp}
\begin{split}
\Delta^2 V_t + k_t^2 V_t & = 0, \\
\Delta^2 V_p + k_p^2 V_p & = 0, \\
\end{split}
\end{equation}
and, by linearity, we take the normal derivative of~\eqref{eq:Vdef} on $\Gamma$ to find boundary conditions
\begin{equation}
\begin{bmatrix}
\tfrac{\partial V_T}{\partial n} \\
\tfrac{\partial V_p}{\partial n}
\end{bmatrix}
 = B
\begin{bmatrix}
g_T \\ g_P
\end{bmatrix} 
=
\begin{bmatrix}
\MCM g_T + t_+(1-i\gamma\Lambda) g_P \\
\MCM g_T + t_-(1-i\gamma\Lambda) g_P
\end{bmatrix}
\equiv
\begin{bmatrix}
g_{V_t} \\ g_{V_p}
\end{bmatrix}.
\end{equation}
So, with $g_T$ and $g_P$ given \emph{a priori}, the scattered field formulation of Morse-Ingard can be solved as a pair of decoupled Helmholtz scattering problems.

\section{Far field boundary conditions}
\label{sec:BC}
\subsection{Boundary conditions for the Helmholtz problem}
We can state far-field boundary conditions for Morse-Ingard and formulate appropriate radiation boundary conditions on $\Sigma$ by transforming such conditions for each of the decoupled Helmholtz equations in~\eqref{eq:scatMIdecomp}.  
So, we begin with the equation
\begin{equation}
\label{eq:helmholtz}
-\Delta u - \kappa^2 u = 0
\end{equation}
posed on $\Omega$, together with Neuman boundary condition
\begin{equation}
\label{eq:helmholtzneumann}
\frac{\partial u}{\partial n} = g
\end{equation}
on the scattering boundary $\Gamma$.
The relevant boundary condition at infinity is the well-known Sommerfeld radiation condition~\cite{colton_inverse_1998,sommerfeld1912greensche}, which requires that
\begin{equation}
\lim_{|x| \rightarrow \infty} |x|^{\frac{n-1}{2}} \left( \frac{\partial}{\partial |x|} - i \kappa \right) u = 0.
\label{eq:sommerfeld-helmholtz}
\end{equation}

A simple approximate boundary condition arises from imposing Sommerfeld at finite radius, i.e.
\begin{equation}
\label{eq:transmissionHelmholtz}
\left( \frac{\partial}{\partial n} - i \kappa \right) u = 0
\end{equation}
on the exterior boundary $\Gamma$.  This is sometimes called the \emph{transmission} boundary condition.  
If the exterior boundary $\Sigma$ is at some radius $R$ away from the origin, this creates an $\mathcal{O}(R^{-2})$ perturbation of $u$ apart from any numerical discretization error, although one may mitigate the computational cost of increasing $R$ by simultaneously increasing the mesh spacing toward the outer boundary~\cite{goldstein1982exterior}.

An alternative approach is the technique of \emph{perfectly matched layers}~\cite{berenger1994perfectly,Druskin2016}, in which one modifies the PDE near the boundary in a so-called \emph{sponge region}.  
The modified coefficients effectively absorb outgoing waves and allow small computational domains, but the resulting algebraic equations do not yield to efficient techniques such as multigrid.
One can use a Schwarz-type method to separately handle the sponge region with a direct solver and the rest of the domain with multigrid or another fast solver~\cite{safin}, although the known method has sub-optimal complexity in three dimensions.

Many \emph{nonlocal} approaches to domain truncation have also been given.  Most classically, the Dirichlet-to-Neumann map (DtN) or Stekhlov-Poincaré operator can be used on the artificial boundary.  By mapping between types of boundary data, DtN operators require, in principle, the solution of a boundary value problem.
In practice, this is often realized by restricting the truncation
boundary to (mappings of) a simple geometry in which separation of variables can be
performed and then making use of a truncated Fourier series.
In~\cite{kirby2021finite}, we give a new approach to nonlocal boundary conditions that requires only the evaluation of non-singular layer potentials, i.e. a surface convolution with the Helmholtz free-space Green's function or its derivatives. A fast algorithm such as the Fast Multipole Method is useful to avoid quadratic complexity. In its continuous (i.e. not yet discretized) form, this boundary condition is exact, and discretization with any order of accuracy is straightforward. In Subsection~\ref{nlbc}, we recall the formulation of this condition for the Helmholtz operator and develop it for the Morse-Ingard system.

\subsection{Far-field conditions for Morse-Ingard}
\label{sec:morse-ingard-farfield}
Each of the decoupled Helmholtz equations in~\eqref{eq:scatMIdecomp} must satisfy the Sommerfeld condition, so that
\begin{equation}
\begin{split}
& \lim_{|x|\rightarrow \infty} |x|^{\frac{n-1}{2}} \left( \frac{\partial}{\partial |x|} - i k_t \right) V_t \\
= & \lim_{|x|\rightarrow \infty} |x|^{\frac{n-1}{2}} \left( \frac{\partial}{\partial |x|} - i k_p \right) V_p = 0.
\end{split}
\end{equation}
In terms of the variables $T$ and $P$, the Morse-Ingard solution must satisfy
\begin{equation}
\label{eq:MISommerfeld}
\begin{split}
& \lim_{|x|\rightarrow \infty} |x|^{\frac{n-1}{2}} \left( \frac{\partial}{\partial |x|} - i k_t \right)\left[ \MCM T + t_{+} (1-i\gamma \Lambda) P \right]\\
= & \lim_{|x|\rightarrow \infty} |x|^{\frac{n-1}{2}} \left( \frac{\partial}{\partial |x|} - i k_p \right) \left[\MCM T + t_- (1-i\gamma \Lambda) P \right] = 0.
\end{split}
\end{equation}
Because of the large imaginary part of $k_t$ the thermal mode attenuates very quickly, so some simple boundary condition could be suitable for $V_t$. 
The acoustic mode, on the other hand, has only a very small imaginary part and so outgoing waves attenuate very slowly.
Since $V_p$ depends on both $T$ and $P$, however, artificial boundary conditions must act on both of these variables. 

We may find an analog to the transmission boundary condition~\eqref{eq:transmissionHelmholtz} by imposing those conditions on each of the decoupled equations, so that we require
\begin{equation}
\label{eq:transdecoupled}
\begin{split}
\left( \frac{\partial}{\partial n} - i k_t \right) V_t & = 0, \\
\left( \frac{\partial}{\partial n} - i k_p \right) V_p & = 0.
\end{split}
\end{equation}
In terms of $T$ and $P$, we have
\begin{equation}
\label{eq:transcoupled}
\begin{split}
\left( \frac{\partial}{\partial n} - i k_t \right) \left( \MCM T + t_+(1-i \gamma \Lambda) P \right) & = 0, \\
\left( \frac{\partial}{\partial n} - i k_p \right) \left( \MCM T + t_-(1-i \gamma \Lambda) P \right) & = 0.
\end{split}
\end{equation}

With some elementary but involved algebraic manipulation, we have
\begin{equation}
\label{eq:transcoupledsolved}
\begin{split}
\MCM \frac{\partial T}{\partial n} & = i \left[ \frac{t_- k_t - t_+ k_p}{t_--t_+} \right] \MCM T
+ i \left[ \frac{t_+ t_- \left( k_t - k_p \right)}{t_--t_+} \right] \left( 1 - i \gamma \Lambda \right) P \\
\left(1 - i \gamma \Lambda\right) \frac{\partial P}{\partial n} & = i \left[ \frac{k_p - k_t}{t_--t_+} \right] \MCM T
+ i \left[ \frac{t_-k_p - t_+ k_t}{t_--t_+} \right] \left( 1 - i \gamma \Lambda \right) P \\
\end{split}
\end{equation}

This boundary condition is of course only an approximation of the actually desired
far-field condition, in the same sense in which \eqref{eq:transmissionHelmholtz} approximates \eqref{eq:sommerfeld-helmholtz}, i.e. it becomes exact as $\Sigma$ moves outward, analogous to the analysis in~\cite{goldstein1982exterior} for Helmholtz.

On the other hand, in prior work~\cite{brennan2015finite, kirby2020optimal}, we had not yet developed the appropriate Sommerfeld condition for Morse-Ingard and used the \emph{ad hoc} boundary conditions
\begin{equation}
\label{eq:wrongbc}
\frac{\partial T}{\partial n} = 0, \ \ \ \frac{\partial P}{\partial n} = i \sqrt{\gamma} P.
\end{equation}
This condition assumes that no heat is transported from the computational domain and an approximate version of \eqref{eq:transmissionHelmholtz} is applied only to the pressure component. Clearly, this boundary condition is quite different from \eqref{eq:transcoupledsolved}.  In particular, the outgoing wave for $V_p$ carries both pressure and temperature with it, thus avoiding reliance on the decay of $T$ for accurate imposition of the boundary condition.

Define
\[
U = \begin{bmatrix} T \\ P \end{bmatrix}
\]
and
\[
B = \begin{bmatrix} \MCM & 0 \\ 0 & 1- i \gamma \Lambda \end{bmatrix}.
\]
Then, boundary conditions~\eqref{eq:transcoupledsolved} and~\eqref{eq:wrongbc} can be written in the form
\begin{equation}
\label{eq:MItransmissiongen}
B \frac{\partial U}{\partial n} = i A U,
\end{equation}
where for~\eqref{eq:transcoupledsolved} we have
\begin{equation}
  A
  = \begin{bmatrix}
    \frac{t_- k_t - t_+ k_p}{t_--t_+} & \left( \frac{t_+ t_- \left( k_t - k_p \right)}{t_--t_+} \right) \left( 1 - i \gamma \Lambda \right) \\
    \frac{k_p - k_t}{t_--t_+}  &
    \left( \frac{t_-k_p - t_+ k_t}{t_--t_+} \right) \left( 1 - i \gamma \Lambda \right)
  \end{bmatrix},
  \label{eq:tcsbc-matrix}
\end{equation}
and for~\eqref{eq:wrongbc},
\begin{equation}
  A = \begin{bmatrix} 0 & 0 \\ 0 & \sqrt{\gamma}\left(1-i\gamma\Lambda\right)
  \end{bmatrix}
  \label{eq:wrongbc-matrix}
\end{equation}


Before proceeding, we offer a brief remark on perfectly matched layers for the Morse-Ingard equations.
In~\cite{safin, safin2018modeling}, it was found that PML must be applied to both the temperature and pressure in order to achieve accurate domain truncation.
Our discussion of transmission boundary conditions shines further light on this observation.
The acoustic mode involves a linear combination of both temperature and pressure, and hence both variables must be damped at the computational boundary in order to avoid spurious reflections.

\subsection{Nonlocal boundary conditions}
\label{nlbc}
Now, we formulate a nonlocal boundary condition based on a Green's integral representation of the solution.
As mentioned, this provides an (in principle) exact boundary condition without the geometric limitations of DtN techniques.
We introduced this technique for the Helmholtz operator in~\cite{kirby2021finite} and now apply it to Morse-Ingard.

For the Helmholtz problem, we let $\mathcal{K}(x, y)$ be the free-space Green's function, which is given by
\begin{equation}
\label{eq:kermit}
  \mathcal K_\kappa(x):=
    \begin{cases}
        \frac{i}{4} H_0^{(1)}(\kappa |x|) & d = 2, \\
        \frac{i}{4\pi|x|} e^{i\kappa |x|} & d = 3.
    \end{cases}
\end{equation}
Here, $H_0^{(1)}$ is the first-kind Hankel function of order 0.
We recall Green's representation theorem~\cite[Thm.~2.5]{colton_inverse_1998} for the Helmholtz equation:
\begin{equation}
u(x) = D_\kappa(u)(x) - S_\kappa(u)(x)\qquad (x\in\Omega),
\end{equation}
where $D_\kappa$ and $S_\kappa$ refer to the double- and single-layer potentials associated with wave number $\kappa$, respectively.  These are
\begin{align}
\label{eq:S}
S_\kappa(u)(x) &= \int_\Gamma \mathcal K_\kappa (x - y ) u(y)dy,\\
\label{eq:D}
D_\kappa(u)(x) &= \operatorname{PV} \int_\Gamma \left( \tfrac{\partial}{\partial n} \mathcal K_\kappa( x-y) \right)u(y) dy.
\end{align}
In anticipation of applying our technique to the decoupled Helmholtz form of Morse-Ingard~\eqref{eq:MIdecomp}, we include the wave number $\kappa$ and use distinct layer potentials with $\kappa=k_p, k_t$.

Using the scattering boundary condition~\eqref{eq:helmholtzneumann} on $\Gamma$, this becomes (omitting the spatial argument)
\begin{equation}
\label{eq:startbc}
u = D_\kappa(u) - S_\kappa(f).
\end{equation}
This representation is valid away from the scattering boundary $\Gamma$, and in particular, on $\Sigma$.  
Hence, we can take its normal derivative, so that on $\Sigma$
\begin{equation}
\label{eq:nonlocal0}
\tfrac{\partial u}{\partial n} = \tfrac{\partial}{\partial n} \left( D_\kappa(u) - S_\kappa(f) \right).
\end{equation}
In~\cite{kirby2021finite}, we subtracted $i \kappa u$ from each side of~\eqref{eq:startbc} and rearranged to arrive at a nonlocal Robin-type boundary condition
\begin{equation}
\label{eq:nonlocal1}
\tfrac{\partial u}{\partial n} =
i \kappa u - \left(i \kappa - \tfrac{\partial}{\partial n} \right) \left( D_\kappa(u) - S_\kappa(f) \right).
\end{equation}
More generally, we can subtract some $i \sigma$ times $u$ from each side of ~\eqref{eq:startbc} to write the condition
\begin{equation}
\label{eq:nonlocalsigma}
\tfrac{\partial u}{\partial n} =
i \sigma u - \left(i \sigma - \tfrac{\partial}{\partial n} \right) \left( D_\kappa(u) - S_\kappa(f) \right), 
\end{equation}
and then~\eqref{eq:nonlocal0} is obtained with $\sigma=0$ and~\eqref{eq:nonlocal1} with $\sigma = \kappa$.

Now, we can apply this boundary condition, in either the form~\eqref{eq:nonlocal0} or~\eqref{eq:nonlocal1} to the decoupled Helmholtz system and back-convert to obtain appropriate nonlocal boundary conditions for Morse-Ingard.  
As in deriving the local transmission condition, we start with the decoupled form and see what is implied in the coupled form.
We let $\boldsymbol{\sigma} = (\sigma_t, \sigma_p)$ be a pair of complex numbers. 
Then, we apply the boundary condition~\eqref{eq:nonlocalsigma} to each equation of~\eqref{eq:MIdecomp} on~$\Sigma$, so that we have
\begin{equation}
\begin{split}
\frac{\partial V_t}{\partial n} & = i \sigma_t V_t - \left( i \sigma_t - \frac{\partial}{\partial n} \right) \left( D_{k_t}(V_t) - S_{k_t}(g_{V_t}) \right),\\
\frac{\partial V_p}{\partial n} & = i \sigma_p V_p - \left( i \sigma_p - \frac{\partial}{\partial n} \right) \left( D_{k_p}(V_p) - S_{k_p}(g_{V_p}) \right).
\end{split}
\end{equation}

Now, we substitute in for $V_t$ and $V_p$ via \eqref{eq:Vdef}, but not in the layer potentials:
\begin{equation}
\label{eq:nonlocbcjustsubbed}
\begin{split}
\tfrac{\partial}{\partial n}\left(\MCM T + t_+(1-i \gamma \Lambda) P\right)
 = & i \sigma_t \left(\MCM T + t_+(1-i \gamma \Lambda) P\right) + R_t, \\
\tfrac{\partial}{\partial n}\left(\MCM T + t_-(1-i \gamma \Lambda) P\right)
 = & i \sigma_p \left(\MCM T + t_-(1-i \gamma \Lambda) P\right) + R_p, \\
\end{split}
\end{equation}
where
\begin{equation}
\begin{split}
R_t & = - \left( i \sigma_t - \tfrac{\partial}{\partial n} \right) D_{k_t}(V_t) + G_t, \\
R_p & = - \left( i \sigma_p - \tfrac{\partial}{\partial n} \right) D_{k_p}(V_p) + G_p,
\end{split}
\end{equation}
and $G_t = \left( i \sigma_t - \tfrac{\partial}{\partial n} \right) S_{k_t}(g_{V_t}) $ and a similar definition for $G_p$.

Similar manipulations leading from~\eqref{eq:transcoupled} to~\eqref{eq:transcoupledsolved} let us rearrange these equations to
\begin{equation}
\label{eq:nonloccoupledsolved}
\begin{split}
\MCM \frac{\partial T}{\partial n} = & i \left[ \frac{t_- \sigma_t - t_+ \sigma_p}{t_--t_+} \right] \MCM T \\
& + i \left[ \frac{t_+ t_- \left( \sigma_t - \sigma_p \right)}{t_--t_+} \right]\left( 1 - i \gamma \Lambda \right) P \\
& + \frac{t_- R_t - t_+ R_p}{t_- - t_+}, \\
\left(1 - i \gamma \Lambda\right) \frac{\partial P}{\partial n} = & i \left[ \frac{\sigma_p - \sigma_t}{t_--t_+} \right] \MCM T \\
& + i \left[ \frac{t_-\sigma_p - t_+ \sigma_t}{t_--t_+} \right] \left( 1 - i \gamma \Lambda \right) P  \\
& + \frac{R_p - R_t}{t_- - t_+}.
\end{split}
\end{equation}
Note that the nonlocality is contained in $R_p$ and $R_t$, each of which depend on both $T$ and $P$ through either $V_t$ or $V_p$.


\section{Variational formulation}
\label{sec:varform}
In this section, we give a finite element formulation of the Morse-Ingard equations~\eqref{eq:scatPDE} under various boundary conditions.  First, we establish some notation.
We let $L^2(\Omega')$ denote the standard space of complex-valued functions with moduli square-integrable over $\Omega'$ and $H^k(\Omega') \subset L^2(\Omega)$ the subspace consisting of functions with weak derivatives up to and including order $k$ also lying in $L^2(\Omega')$.  For any Banach space $\mathcal{V}$, we let $\| \cdot \|_{\mathcal{V}}$ denote its norm, with the subscript typically omitted when $\mathcal{V} = L^2(\Omega')$.  

The space $L^2(\Omega')$ is equipped with the standard inner product
\begin{equation}
\label{eq:l2ip}
(f, g) = \int_{\Omega'} f(x) \overline{g(x)} dx,
\end{equation}
and we also define the inner product over a portion of the boundary $\tilde{\Gamma} \subset \partial \Omega'$ by
\begin{equation}
\langle f, g \rangle_{\tilde{\Gamma}} = \int_{\tilde{\Gamma}} f(s) \overline{g(s)} ds.
\end{equation}

We partition $\Omega'$ into a family of conforming, quasi-uniform triangulations~\cite{BrennerScott} $\{ \mathcal{T}_h \}_{h > 0}$.  
Let $\mathcal{V}_h$ be the standard space of continuous piecewise polynomials of some degree $k\geq 1$ over $\mathcal{T}_h$.  
Since we are dealing with a system of two PDE, we define $\boldsymbol{\mathcal{V}}_h = \mathcal{V}_h \times \mathcal{V}_h$.

We multiply the equations of~\eqref{eq:scatPDE} by the conjugate of the test functions $v, w \in \mathcal{V}_h$, respectively and integrate by parts over $\Omega'$.  We let $U=(T, P)$ and $\Psi=(v,w)$.  Applying the Neumann boundary conditions~\eqref{eq:scatBC} on $\Gamma$ but not taking action yet on $\Sigma$, we have
\begin{equation}
  \label{eq:preweakform}
  \begin{split}
    \MCM \left( \nabla T, \nabla v \right) 
    - \langle \MCM \tfrac{\partial T}{\partial n}, v \rangle_{\Sigma} - i \left( T, v \right)
    + i \tfrac{\gamma-1}{\gamma} \left( P, v \right) & = 
    \langle g_T, v \rangle_{\Gamma}, \\
   \gamma\left( 1 - \tfrac{\Lambda}{\MCM} \right) \left( T, w \right)
   + \left( 1 - i \gamma \Lambda \right)
   \left[ \left(\nabla P, \nabla w\right)
     - \langle \tfrac{\partial P}{\partial n} , w \rangle_{\Sigma} \right] \\
  - \left[ \gamma \left(1-\tfrac{\Lambda}{\MCM}\right)+
     \tfrac{\Lambda}{\MCM}\right]\left(P, w\right)
  &  = \langle g_P , w \rangle_{\Gamma}.
  \end{split}
\end{equation}
We add these equations together and define $a_0: \mathbf{V} \times \mathbf{V} \rightarrow \mathbb{C}$ as consisting of the volumetric terms:
\begin{equation}
\begin{split}
a_0\left( U, \Psi \right)
= &     \MCM \left( \nabla T, \nabla v \right)  - i \left( T, v \right) + i \tfrac{\gamma-1}{\gamma} \left( P, v \right) \\
& +   \gamma\left( 1 - \tfrac{\Lambda}{\MCM} \right) \left( T, w \right)
   + \left( 1 - i \gamma \Lambda \right)
   \left(\nabla P, \nabla w\right) \\
 & - \left[ \gamma \left(1-\tfrac{\Lambda}{\MCM}\right)+
     \tfrac{\Lambda}{\MCM}\right]\left(P, w\right),
\end{split}
\end{equation}
and $F_0$ involving those boundary terms on the right-hand side:
\begin{equation}
F_0(\Psi) =
\langle g_T, v \rangle_{\Gamma} + \langle g_P, w \rangle_{\Gamma}.
\end{equation}
Then, we can write~\eqref{eq:preweakform} as
\begin{equation}
\label{eq:preweakform2}
a_0(U, \Psi) - \langle \MCM \tfrac{\partial T}{\partial n}, v \rangle_{\Sigma}
- (1-i \gamma \Lambda) \langle \tfrac{\partial P}{\partial n}, w \rangle_{\Sigma}
= F_0(\Psi).
\end{equation}

At this point, we can close the system by selecting any of the boundary conditions discussed in Section~\ref{sec:BC} and substituting in the relevant expressions for $\tfrac{\partial T}{\partial n}$ and $\tfrac{\partial P}{\partial n}$.  Following the general form of the local boundary condition~\eqref{eq:MItransmissiongen}, we define $a_A$ by
\begin{equation}
\begin{split}
a_A(U, \Psi) = & a_0(U, \Psi)  - \langle \alpha_{11} T + \alpha_{12} P, v \rangle_{\Sigma} 
\\ & -  \langle  \alpha_{21} T + \alpha_{22} P, w \rangle_{\Sigma}
\end{split}
\end{equation}
with
\[
A=\begin{bmatrix}
  \alpha_{11} & \alpha_{12}\\
  \alpha_{21} & \alpha_{22}
\end{bmatrix},
\]
for the respective $A$ chosen, cf. \eqref{eq:tcsbc-matrix} and \eqref{eq:wrongbc-matrix}.
This leads to the variational problem of finding $U \in \mathbf{V}$ such that
\begin{equation}
a_A(U, \Psi)
= F_0(\Psi)
\end{equation}
for all $\Psi \in \mathbf{V}$.

We now consider variational problems corresponding to the exact, but nonlocal, boundary conditions.  
Recall that these boundary conditions are parametrized over the choice of $\boldsymbol{\sigma} = (\sigma_t, \sigma_p)$.
Using boundary condition~\eqref{eq:nonloccoupledsolved} in~\eqref{eq:preweakform} motivates defining the bilinear form
\begin{equation}
  \label{eq:asig}
\begin{split}
a_{\boldsymbol{\sigma}}(U, \Psi) 
 = & a_0(U, \Psi)  \\
& +  \langle \tfrac{1}{t_--t_+} \left[ t_- \left( i \sigma_t - \tfrac{\partial}{\partial n} \right)  D_{k_t}(V_t) - t_+ \left( i \sigma_p - \tfrac{\partial}{\partial n} \right) D_{k_p}(V_p)\right], v \rangle_{\Sigma} \\
& + \langle \tfrac{1}{t_--t_+} \left( i \sigma_p - \tfrac{\partial}{\partial n} \right) D_{k_p}(V_p)
  - \tfrac{1}{t_--t_+} \left( i \sigma_p - \tfrac{\partial}{\partial n} \right) D_{k_t}(V_t), w \rangle_\Sigma, \\
  \equiv & a_0(U, \Psi) + a^{NL}_{\boldsymbol{\sigma}}(U, \Psi)
\end{split}
\end{equation}
and linear form
\begin{equation}
F_{\boldsymbol{\sigma}}(\Psi) = F_0(\Psi)
+ \langle \tfrac{t_-G_t - t_+G_p}{t_--t_+} , v \rangle_{\Sigma} 
- \langle \tfrac{G_t-G_p}{t_--t_+}, w \rangle_{\Sigma}.
\end{equation}

Then, we pose the variational problem of finding $U^{\boldsymbol{\sigma}} \in \boldsymbol{\mathcal{V}}$ such that 
\begin{equation}
a_{\boldsymbol{\sigma}}(U^{\boldsymbol{\sigma}}, \Psi) = F_{\boldsymbol{\sigma}}(\Psi)
\end{equation}
for all $\Psi \in \boldsymbol{\mathcal{V}}$.
In fact, for any choice of $\boldsymbol{\sigma}$, $U^{\boldsymbol{\sigma}}$ is the solution to~\eqref{eq:scatPDE} with scattering boundary conditions~\eqref{eq:scatBC} and the Sommerfeld-type far field condition~\eqref{eq:MISommerfeld}.

A standard Galerkin discretization of this problem is obtained by restricting the test function $\Psi$ to $\boldsymbol{\mathcal{V}}_h$, seeking $U^{\boldsymbol{\sigma}}_h \in \boldsymbol{\mathcal{V}}_h$ such that
\begin{equation}
\label{eq:nonlocdiscretevar}
a_{\boldsymbol{\sigma}}(U^{\boldsymbol{\sigma}}_h, \Psi_h) = F_{\boldsymbol{\sigma}}(\Psi_h)
\end{equation}
for all $\Psi_h \in \boldsymbol{\mathcal{V}}_h$.

Analyzing this discretization follows along the lines proposed in~\cite{kirby2021finite} for the scalar Helmholtz problem -- one establishes a G{\aa}rding inequality for the variational form, by which existing theory~\cite{BrennerScott} for Galerkin methods for elliptic operators provides solvability and optimal $H^1$ and $L^2$ and error estimates, subject to a sufficiently fine mesh.  
We have proven a G{\aa}rding-type inequality for the local form of Morse-Ingard in~\cite{kirby2020optimal}, and the same techniques used to handle the nonlocal terms for Helmholtz in~\cite{kirby2021finite} can be used for Morse-Ingard.  Consequently,
\begin{theorem}
There exists some $h_0 > 0$ such that for $h \leq h_0$, the variational problem~\eqref{eq:nonlocdiscretevar} has a unique solution $U_h^{\boldsymbol{\sigma}}$, and this solution satisfies the best approximation result
\begin{equation}
\left\| U^{\boldsymbol{\sigma}} - U^{\boldsymbol{\sigma}}_h \right\|_{(H^1(\Omega'))^2} \leq C 
\inf_{W_h \in \boldsymbol{\mathcal{V}}_h}
\left\| U^{\boldsymbol{\sigma}} - W_h \right\|_{(H^1(\Omega'))^2},
\end{equation}
and $L^2$ error estimate
\begin{equation}
\left\| U^{\boldsymbol{\sigma}} - U^{\boldsymbol{\sigma}}_h \right\|_{(L^2(\Omega'))^2} \leq C h
\left\| U^{\boldsymbol{\sigma}} -  U^{\boldsymbol{\sigma}}_h \right\|_{(H^1(\Omega'))^2},
\end{equation}
where the constants $C$ in the two inequalities differ from each other but are independent of $h$.
\end{theorem}

\section{Linear algebra}
\label{sec:precond}
A major feature of our nonlocal boundary condition is the opportunity for efficient solvers.  For the Helmholtz problem in~\cite{kirby2021finite}, we demonstrated empirically that preconditioning the entire operator with the local part led to very low GMRES iteration counts.  This, of course, means the cost of inverting the local part of the operator drives the overall cost.  It is well known that high wave numbers lead to notorious difficulty for iterative methods~\cite{ernst2012difficult}.
Fortunately, this is not the case for the parameter regime of interest for Morse-Ingard.  In decoupled form~\eqref{eq:MIdecomp}, the thermal mode has a large imaginary part to complement the large real part,  while the thermal mode has wave number approximately 1.  
Standard multigrid algorithms handle both of these situations effectively~\cite{gander2015applying}, so solving the problem in decoupled form proceeds along lines given in~\cite{kirby2021finite}, followed by forming $T$ and $P$ from $V_t$ and $V_p$.

For the fully coupled system, we introduce a basis $\left\{ \psi_i \right\}_{i=1}^{\dim \boldsymbol{\mathcal{V}}_h}$, and then we can write~\eqref{eq:nonlocdiscretevar} as a linear system
\begin{equation}
A \mathbf{x} = \mathbf{b},
\end{equation}
where
\begin{equation}
\begin{split}
A_{ij} & = a_{\boldsymbol{\sigma}}(\psi_j, \psi_i), \\
\mathbf{b} & = F_{\boldsymbol{\sigma}}(\psi_j),
\end{split}
\end{equation}
and then $U^{\boldsymbol{\sigma}} = \sum_{i=1}^{\dim \boldsymbol{\mathcal{V}}_h} \mathbf{x}_i \psi_i$.

Following the partition of $a_{\boldsymbol{\sigma}}$ given in~\eqref{eq:asig}, we can write
$A = A^L + A^{NL}$, where
\begin{equation}
\begin{split}
A^L_{ij} & = a_0(\psi_j, \psi_i), \\
A^{NL}_{ij} & = a^{NL}_{\boldsymbol{\sigma}}(\psi_j, \psi_i),
\end{split}
\end{equation}
corresponding to the local and nonlocal parts of the bilinear form.  

We can effectively solve the linear using preconditioned GMRES~\cite{saad1986gmres}, which is a parameter-free algorithm approximating the solution of a linear system $A \bfx = \bfb$ as the element of the Krylov subspace $\operatorname{span}\{ A^i \bfb \}_{i=0}^m$ minimizing the equation residual.  
Building the subspace does not require access to entries of $A$, just the action of $A$ on vectors.  
Unlike conjugate gradients, GMRES is not restricted to operators that are symmetric and positive definite.

For most problems arising in the discretization of PDE, GMRES is most frequently used in conjunction with a \emph{preconditioner}.  Mathematically, we multiply the linear system through by some matrix $\widehat{P}^{-1}$:
\begin{equation}
\widehat{P}^{-1} A \bfx = \widehat{P}^{-1} \bfb,
\end{equation}
and so the Krylov space then is $\operatorname{span}\{ \left( \widehat{P}^{-1}A \right)^i \widehat{P}^{-1} \bfb \}_{i=0}^m$.

The overall performance of GMRES typically depends on two factors -- the cost of building and applying the operators $\widehat{P}^{-1}$ and $A$, and the total number of iterations.  
One hopes to obtain a per-application cost that scales linearly (or log-linearly) with respect to the number of unknowns in the linear system, and a total number of GMRES iterations that is bounded independently of the number of unknowns.  We think of $\widehat{P}^{-1}$ being an approximation to the inverse of some matrix $P$ that approximates $A$.
As with the Helmholtz problem in~\cite{kirby2021finite}, we will take $P = A^L$, the local part of the operator.  Then, applying $\widehat{P}^{-1}$ might correspond to applying the inverse of $P$ by a sparse direct method, an application of some block preconditioner such as that analyzed in~\cite{kirby2020optimal}, or some other strategy.

As a partial justification of our choice of preconditioning matrix, when $\widehat{P} = A^L$ so that the inverse is applied exactly, we arrive at a preconditioned matrix of the form
\begin{equation}
\widehat{P}^{-1} A
= \left( A^L \right)^{-1} \left( A^L + A^{NL} \right)
= I + \left( A^L \right)^{-1} A^{NL}.
\end{equation}
Because $A^{NL}$ discretizes a compact operator (layer potential in weak form) and, moreover, $\left(A^L\right)^{-1}$ discretizes the inverse of an elliptic operator, the preconditioned matrix has the form of a (discretization of) a compact perturbation of the identity.  
We suggested obtaining such a form via preconditioning as a heuristic in~\cite{Howle}.  Moret~\cite{moret1997note} gives rigorous GMRES convergence estimates for this situation once one establishes certain bounds on the operators.  
In practice, one might not apply the inverse of $A^L$ exactly.  
Supposing it is spectrally equivalent, one still obtains a compact perturbation of a bounded operator.  Recently, Blechta~\cite{blechta2021stability} has extended Moret's results to describe GMRES convergence for this setting as well.

\section{Numerical results}
\label{sec:numres}
In this section, we apply our various formulations of the Morse-Ingard equations.
All of our numerical experiments are conducted using Firedrake~\cite{Rathgeber:2016}.  We generate meshes with Gmsh~\cite{geuzaine2009gmsh}, using its internal interface to OpenCascade for constructive solid geometry. 
Our Firedrake experiments make use of the recently-developed \lstinline{ExternalOperator} capability, which is a foreign function interface allowing users to incorporate operators implemented by other packages into UFL expressions~\cite{Alnaes:2014}.  This technology will be documented elsewhere.

We use this to enable Firedrake to express layer potentials and  evaluate them using the \pytential{} package.
\pytential{}~\cite{pytential} is an open-source, MIT licensed software system for evaluating layer potentials from source geometry represented by
unstructured meshes with high accuracy and near-optimal complexity.
\pytential{} provides for the discretization of a source surface using tools for high-order accurate nonsingular quadrature~\cite{xiao_numerical_2010, meshmode}, its refinement according to accuracy requirements~\cite{wala_fast_2019},
and, finally, the evaluation of integral operators via quadrature by expansion
(QBX)~\cite{kloeckner_quadrature_2013} and the associated GIGAQBX fast
algorithm~\cite{wala_fast_2018}, with rigorous accuracy guarantees in
two and three dimensions~\cite{wala_approximation_2020}. This fast
algorithm, can, in turn make use of \software{FMMLIB}~\cite{gimbutas_fast_2009,pyfmmlib}
for the evaluation of translation operators in the moderate-frequency
regime for the Helmholtz operator.  In our two-dimensional experiments, we use an FMM order of 15, which provides sufficient accuracy for the accuracy of layer potential evaluation to not limit the overall accuracy obtained.
While the integrals in our variational problem do not require the singular integral technology provided in QBX, it does provide robustness in the case of $\Sigma$ and $\Gamma$ are chosen to lie close together.

Our simulations are performed on a Debian Linux machine using dual Broadwell-EP 12-core processors running at 2.20GHz.  The machine has 256GB of 2133GHz DDR4 RAM.   In these simulations, \firedrake{} and \petsc{} were using a single core,  \pytential{} evaluates the layer potentials in multicore mode through the PoCL implementation of OpenCL.

\subsection*{Accuracy of the decoupled formulation}
\begin{figure}[ht]
\begin{subfigure}[l]{0.5\textwidth}
\centering
\includegraphics[width=\textwidth]{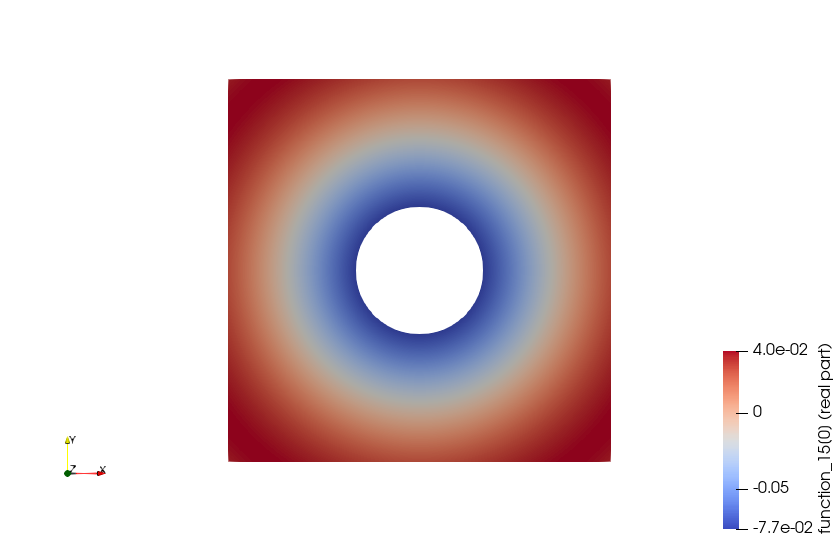}
\caption{Real part of $G_T$}
\end{subfigure} \hfill
\begin{subfigure}[l]{0.5\textwidth}
\centering
\includegraphics[width=\textwidth]{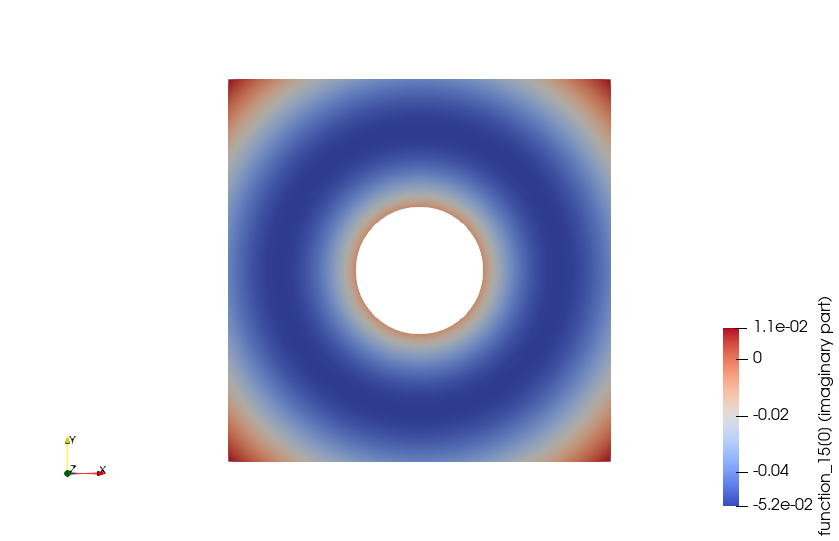}
\caption{Imaginary part of $G_T$}
\end{subfigure} \\
\begin{subfigure}[l]{0.5\textwidth}
\centering
\includegraphics[width=\textwidth]{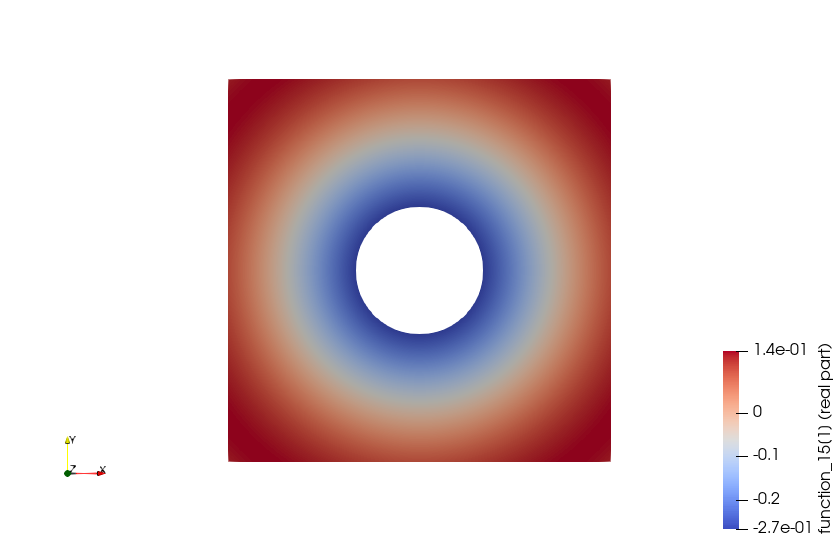}
\caption{Real part of $G_P$}
\end{subfigure} \hfill
\begin{subfigure}[l]{0.5\textwidth}
\centering
\includegraphics[width=\textwidth]{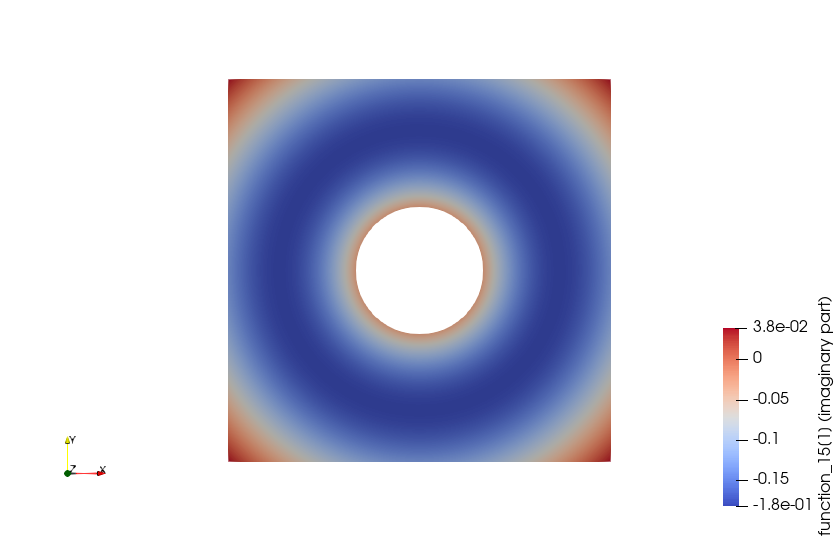}
\caption{Imaginary part of $G_P$}
\end{subfigure}
\caption{Real and imaginary parts of the temperature/pressure Green's functions~\eqref{eq:kermit}.}
\label{fig:plotg}
\end{figure}

\begin{figure}[ht]
\begin{subfigure}[l]{0.5\textwidth}
\centering
\includegraphics[width=\textwidth]{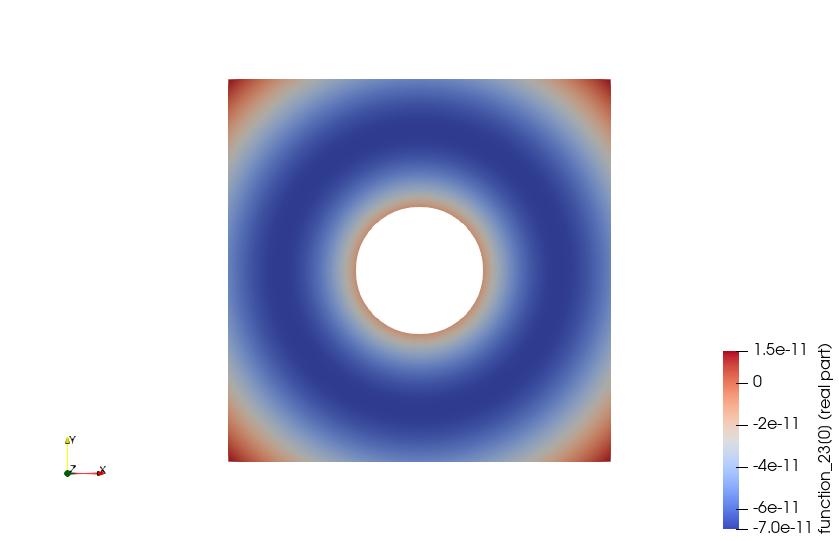}
\caption{Real part of $G_{V_t}$}
\end{subfigure} \hfill
\begin{subfigure}[l]{0.5\textwidth}
\centering
\includegraphics[width=\textwidth]{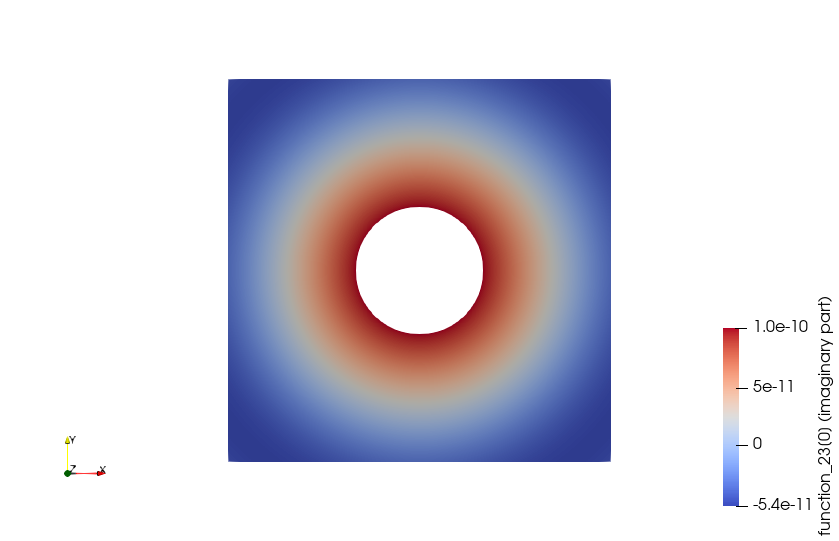}
\caption{Imaginary part of $G_{V_t}$}
\end{subfigure} \\
\begin{subfigure}[l]{0.5\textwidth}
\centering
\includegraphics[width=\textwidth]{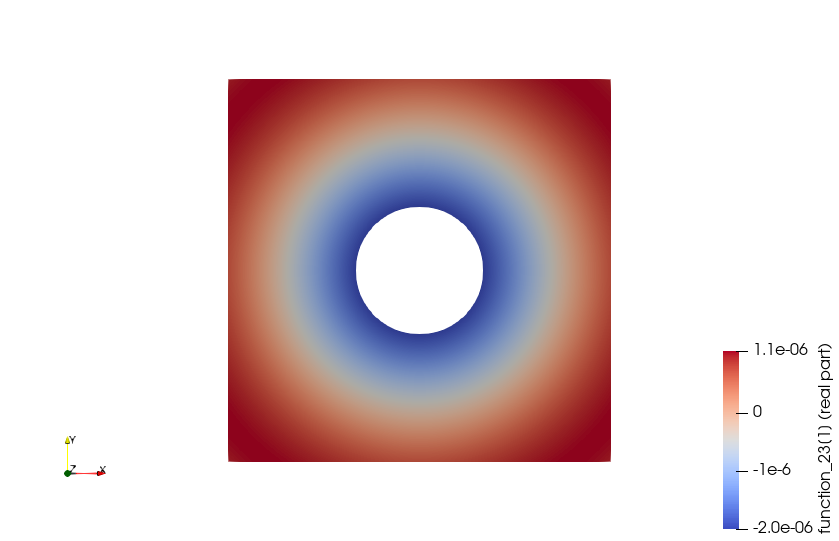}
\caption{Real part of $G_{V_p}$}
\end{subfigure} \hfill
\begin{subfigure}[l]{0.5\textwidth}
\centering
\includegraphics[width=\textwidth]{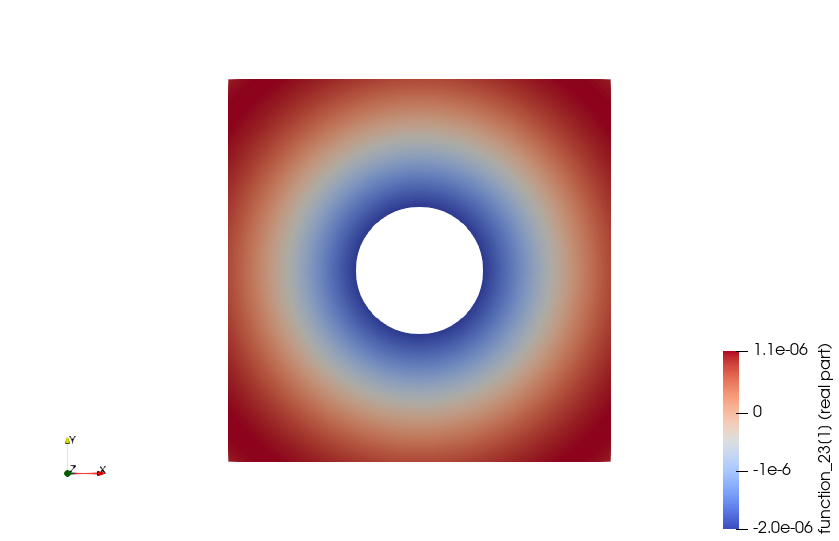}
\caption{Imaginary part of $G_{V_p}$}
\end{subfigure}
\caption{Real and imaginary parts of the thermal and acoustic modes corresponding to the Green's functions in Figure~\ref{fig:plotg}.}
\label{fig:plotvg}
\end{figure}

First, we demonstrate the relative accuracy obtainable through the coupled and decoupled forms of the Morse-Ingard equations.  In these cases, we consider a simple square domain $[-1.5,1.5]^2$ with a circle of radius 2/3 removed from the center.  We selected Neumann boundary data on $\Gamma \cup \Sigma$ as $g_T$ and $g_P$ being the normal derivatives of the temperature and pressure Green's functions centered at the origin.  We plot the real and imaginary parts of these Green's functions in Figure~\ref{fig:plotg}.  For comparison, we also plot the equivalent thermal and acoustic modes $V_t$ and $V_p$ in Figure~\ref{fig:plotvg}.  These are the solutions to the pair of decoupled Helmholtz equations~\eqref{eq:scatMIdecomp}.

\begin{figure}[ht]
\centering
\begin{tikzpicture}[scale=0.88]
\begin{loglogaxis}[xlabel={$N_v$}, ylabel={$L^2$ error}]
\addplot[dashed, red]
table [x=NV,y=E, col sep=comma]{decoupled.dat};
\addlegendentry{Decoupled}
\addplot[dotted, blue]
table [x=NV,y=E, col sep=comma]{coupled.dat};
\addlegendentry{Coupled}
\addplot[dashdotted, green]
table [x=NV,y=E, col sep=comma]{single_eq.dat};
\addlegendentry{Single eq}

\end{loglogaxis}
\end{tikzpicture}
\caption{Relative $L^2$ accuracy of solving Morse-Ingard equations with Neumann boundary conditions in both coupled and decoupled forms using quadratic Lagrange elements, and in approximating the solution by only solving for the thermal mode.  On coarse meshes, the solutions for the two systems agree, but the convergence in the decoupled form levels off.  The approximation using only the thermal mode is slightly worse than the systems on coarse meshes but achieves comparable accuracy to the decoupled form on finer meshes.}
\label{fig:coupledvsdecoupled}
\end{figure}
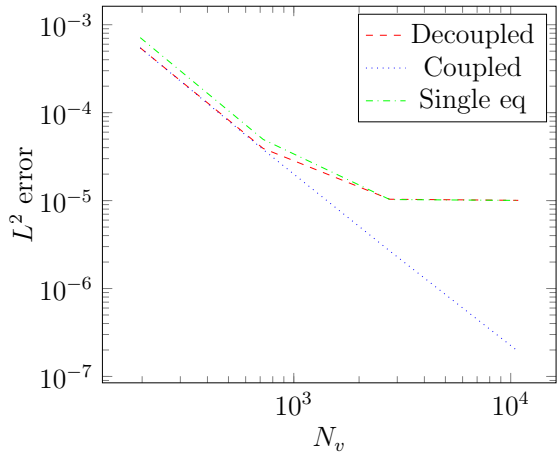

For this experiment, we consider a mesh of quadratic triangles and discretize the Morse-Ingard equations in decoupled form~\eqref{eq:scatMIdecomp} using quadratic Lagrange elements for both $V_t$ and $V_p$.  After solving the decoupled equations, we compute the equivalent $T$ and $P$.  We plot the relative $L^2$ error in the resulting solution in Figure~\eqref{fig:coupledvsdecoupled}

If we repeat the experiment with cubic approximation to $V_t$ and $V_p$, we see similar results, with the maximal accuracy of about $10^{-5}$ reached on an even coarser mesh.  Similary, we can use linear discretizations and observe second order convergence until this threshold is reached on much finer meshes.  This suggests an intrinsic limit to accuracy obtainable in the decoupled form, resulting either from the conditioning of the transformation~\eqref{eq:Vdef} or some difficulty in accurately solving for the thermal mode.
For comparison, we plot the error obtained in solving the same problem in the coupled form of the equations along with the error in the decoupled form in Figure~\ref{fig:coupledvsdecoupled}.

Although the barrier to accuracy in the decoupled form is disappointing, there is some compensation.  For physical parameters of interest, the thermal wave attenuates extremely quickly away from the boundary. 
Moreover, by comparing the top row of Figure~\eqref{fig:plotvg} to the bottom, we see that the thermal mode is roughly 4-5 orders of magnitude smaller than the acoustic mode.  From this, it seems that 4-5 digits of relative accuracy could be obtained simply by solving~\eqref{eq:scatMIdecomp} for the acoustic mode $V_p$, approximating $V_t \approx 0$, and then constructing the pressure and temperature.
This simply requires solving a wave equation with wave number essentially 1 (plus a very small imaginary part).  The error in making this approximation is also shown Figure~\eqref{fig:coupledvsdecoupled}

\subsection*{Accuracy of boundary conditions}
\begin{figure}
\centering
\begin{tikzpicture}[scale=10]
\draw (-0.1125, 0) rectangle (0.1125, 0.723);
\draw[fill=yellow!30] (-0.015, 0.673) -- (-0.075, 0.673) --
(-0.075, 5e-2) -- (0.075, 5e-2) -- (0.075, 0.673) -- (0.015, 0.673) -- (0.015, 0.298)
arc[start angle=0, end angle=-180, radius=0.015] -- cycle;
\node at (0.135, 0.35) {$\Sigma$};
\node at (0.0, 0.025) {$\Gamma$};
\draw[fill=red] (-0.0375, 0.1665) circle (0.005);
\end{tikzpicture}
\caption{Computational tuning fork domain is the rectangle minus yellow shaded region.  $\Sigma$ is the outer rectangle, and $\Gamma$ is the boundary of the tuning fork itself.  The red circle shows the location of the point source.}
\label{fig:tfsource}
\end{figure}
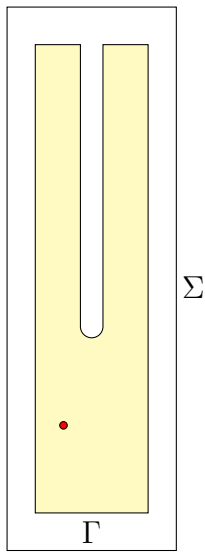

Now, for assessing the accuracy attainable for various boundary condition formulations, we turn to a more realistic tuning fork geometry, showin in Figure~\ref{fig:tf}.  As a manufactured solution, we choose the free space Green's function associated with a point source outside the computational domain.  In our case, we choose the point located midway up the tuning fork base, and off of the vertical line of symmetry, as shown by the red dot in Figure~\ref{fig:tfsource}.

We consider three different kinds of boundary conditions -- the ``wrong'' transmission condition~\eqref{eq:wrongbc}, the ``right'' transmission boundary condition~\eqref{eq:transcoupledsolved}, and the nonlocal boundary conditions~\eqref{eq:nonloccoupledsolved}.  Just as for the Helmholtz problem, the local boundary conditions represent a fundamental perturbation of the boundary value problem, and the finite element solution cannot converge to the correct solution.  The solutions to the problem with~\eqref{eq:wrongbc} converge to a wrong answer with about 50.9\% relative error.  Using~\eqref{eq:transcoupledsolved} is only slightly less bad, as numerical converge to a solution with about 45.5\% relative error.  As seen in Figure~\ref{fig:tf}, we have a relatively small domain enclosing the tuning fork so that the effects of reflected waves at the boundary is considerable.  Following Goldstein's technique for Helmholtz~\cite{goldstein1982exterior} one could obtain increased accuracy by increasing the domain size and increasing element size near the outer boundary~$\Sigma$, but the $\mathcal{O}(R^{-2})$ convergence rate he proves means the domain would need to be considerably larger to achieve a convergent method.

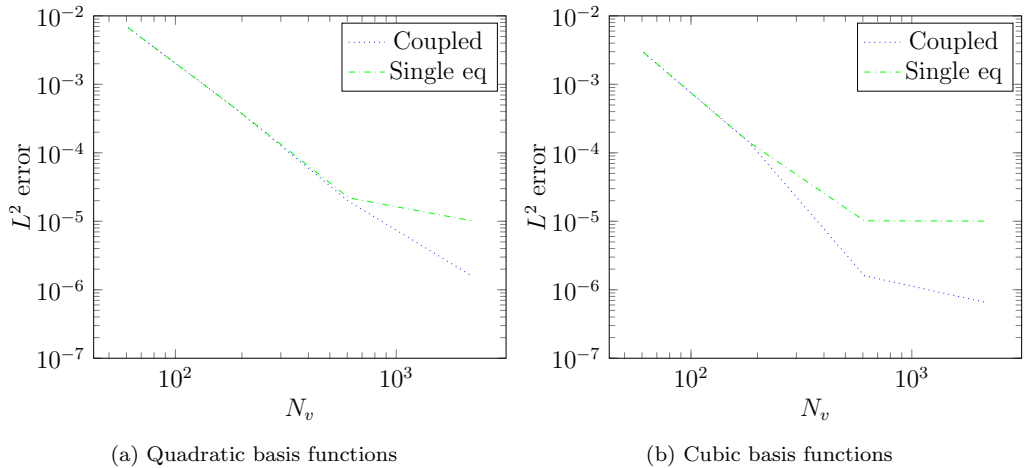
\begin{figure}[ht]
  \centering
  \begin{subfigure}[l]{0.49\textwidth}
    \begin{tikzpicture}[scale=0.8]
      \begin{loglogaxis}[xlabel={$N_v$}, ylabel={$L^2$ error},
          ymax=0.01, ymin=1.e-7]
        \addplot[dotted, blue]
        table [x=NV,y=E, col sep=comma]{coupled_nlbc2.dat};
        \addlegendentry{Coupled}
        \addplot[dashdotted, green]
        table [x=NV,y=E, col sep=comma]{single_eq_nlbc2.dat};
        \addlegendentry{Single eq}
      \end{loglogaxis}
    \end{tikzpicture}
    \caption{Quadratic basis functions}
  \end{subfigure}
  \begin{subfigure}[l]{0.49\textwidth}
    \begin{tikzpicture}[scale=0.8]
      \begin{loglogaxis}[xlabel={$N_v$}, ylabel={$L^2$ error},
          ymax=0.01, ymin=1.e-7]
        \addplot[dotted, blue]
        table [x=NV,y=E, col sep=comma]{coupled_nlbc3.dat};
        \addlegendentry{Coupled}
        \addplot[dashdotted, green]
        table [x=NV,y=E, col sep=comma]{single_eq_nlbc3.dat};
        \addlegendentry{Single eq}
        
      \end{loglogaxis}
    \end{tikzpicture}
    \caption{Cubic basis functions}
  \end{subfigure}  
\caption{Relative $L^2$ accuracy of solving Morse-Ingard equations with nonlocal  boundary conditions.  We compare the fully coupled model~\eqref{eq:scatPDE} to the solution obtained by solving only for the acoustic mode $V_p$ and approximating the thermal mode $V_t \approx 0$.}
\label{fig:nlbcacc}
\end{figure}

On the other hand, we can use the small computational domain in Figure~\ref{fig:tf} in conjunction with nonlocal boundary conditions to achieve a very accurate  solution, as shown in Figure~\ref{fig:nlbcacc}.  We see comparable accuracy to that obtained in Figure~\ref{fig:coupledvsdecoupled}.  We note that on coarse meshes, solving the coupled system and the approximating the thermal mode with 0 give comparable solutions, although with finer meshes the error in this approximation becomes more apparent.

Because our method gives quite high accuracy on rather coarse meshes, with but a few thousand vertices, we have not pursued iterative methods for the local part of the Morse-Ingard operator such as those given in~\cite{kirby2020optimal}.
Instead, we perform a sparse factorization of the local part of the operator $A^L$ and use this as a preconditioner for the full matrix $A$.
For larger problems, especially in three dimensions, this will become impractical.
Still, in all of our simulations, we observed between 13 and 15 GMRES iterations to reach a relative Euclidean norm tolerance of $10^{-12}$.
This provides a lower bound of what one might expect of a nested iteration or simply using a preconditioner for $A^L$ rather than $(A^L)^{-1}$ as a preconditioner for $A$.



\section{Conclusions}
\label{sec:conc}
We have developed exact truncating boundary conditions for the Morse-Ingard equations.
These boundary conditions use a Green's formula representation of the solution in terms of layer potentials and work in general unstructured geometry.
The action of the discrete operators may be evaluated efficiently using matrix-free finite elements and a fast multipole method for the layer potentials, and the linear system may be effectively preconditioned with the local part of the operator.
Standard convergence theory holds for the Galerkin discretization, and the method gives good accuracy on small computational domains even with relatively coarse meshes.

In the future, we hope to pursue a rigorous suite of three-dimensional calculations, compute with iterative treatment for $A^L$, especially as ongoing \pytential{} improves its performance for three-dimensional problems.
We also hope to study models in which the Morse-Ingard are equations are coupled to the tuning fork displacement.

\bibliographystyle{elsarticle-num} 
\bibliography{myBib}





\end{document}